\documentclass[11pt]{article}

\usepackage{authblk}
\usepackage{amsfonts}
\usepackage{amsthm}
\usepackage{amsmath}
\usepackage{graphicx}
\usepackage[usenames,dvipsnames]{xcolor}
\usepackage{enumerate}
\usepackage{caption}
\usepackage{doi}

\newcommand{\DD}{\mathbb{D}}
\newcommand{\R}{\Bbb{R}}
\newcommand{\cadlag}{c\`adl\`ag}
\newcommand{\dto}{\downarrow}
\newcommand{\toj}{\stackrel{\mathrm{J_1}}{\longrightarrow}}

\newcommand{\myAA}{\mathbb{A}}
\newcommand{\PP}{\textbf{P}}
\newcommand{\EE}{\textbf{E}}
\newcommand{\wt}{\widetilde}
\newcommand{\halmos}{\quad\hfill\mbox{$\Box$}}
\newcommand{\sko}{\stackrel{J_1}{\longrightarrow}}
\newcommand{\pibar}{\overline{\Pi}}
\newcommand{\pibarinv}{\overline{\Pi}^{\leftarrow}}
\newcommand{\pibarpinv}{\overline{\Pi}^{+,\leftarrow}}
\newcommand{\todr}{\stackrel{\mathrm{D}}{\longrightarrow}}
\newcommand{\rmd}{{\rm d}}

\newcommand{\eqdr}{\stackrel{\mathrm{D}}{=}}

\newcommand{\BB}{\mathbb{B}}
\newcommand{\NN}{\mathbb{N}}
\newcommand{\RR}{\mathbb{R}}

\newcommand{\EEE}{{\cal E}}

\newcommand{\CCC}{{\cal C}}

\newcommand{\FFF}{{\cal F}}

\newcommand{\VVV}{{\cal V}}
\newcommand{\SSS}{{\cal S}}
\newcommand{\TTT}{{\cal T}}

\newcommand{\YYY}{{\cal Y}}
\newcommand{\III}{{\cal I}}

\newcommand{\RRR}{{\cal R}}

\newcommand{\TTTT}{\mathfrak{T}}
\newcommand{\CCCC}{\mathfrak{C}}

\newcommand{\be}{\begin{equation}}
\newcommand{\ee}{\end{equation}}
\newcommand{\bea}{\begin{eqnarray}}
\newcommand{\eea}{\end{eqnarray}}
\newcommand{\bean}{\begin{eqnarray*}}
\newcommand{\eean}{\end{eqnarray*}}
\newcommand{\ben}{\begin{equation*}}
\newcommand{\een}{\end{equation*}}
\newcommand{\ba}{\begin{aligned}}
\newcommand{\ea}{\end{aligned}}

\newtheorem{theorem}{Theorem}[section]
 
\newtheorem{prop}{Proposition}[section]
\newtheorem{lemma}{Lemma}[section]
\newtheorem{remark}{Remark}[section]
\newtheorem{example}{Example}[section]

\def\boxit#1{\vbox{\hrule\hbox{\vrule\kern6pt
          \vbox{\kern6pt#1\kern6pt}\kern6pt\vrule}\hrule}}

\newcommand*{\medcup}{\mathbin{\scalebox{0.95}{\ensuremath{\bigcup}}}}%

\title{Functional Laws for Trimmed L\'evy Processes} 

\author{Boris Buchmann\thanks{boris.buchmann@anu.edu.au}}
\author{Yuguang F. Ipsen\thanks{Corresponding author: yuguang.fan@anu.edu.au}}
\author{Ross A. Maller\thanks{ross.maller@anu.edu.au}}
\affil{Research School of Finance,  Actuarial Studies \&  Statistics, 
the Australian National University, 26C Kingsley St, Acton ACT 2601, Australia.}

\begin{document}

\maketitle

\begin{abstract}
Two different ways of trimming the sample path of a stochastic process in  $\DD[0,1]$:
global (``trim as you go") trimming and record time (``lookback") trimming  are analysed to find conditions for the corresponding operators to be continuous with respect to the (strong) $J_1$-topology.
A key condition is that there should be no ties among the largest ordered jumps of the limit process. As an application of the theory, via the continuous mapping theorem we prove limit theorems for trimmed L\'evy processes, using the functional convergence of the underlying process to a stable process. The results are applied to a reinsurance ruin time problem. 

\end{abstract}


\section{Introduction}\label{s0}
By ``trimming" a process we mean identifying ``large" jumps of the process, in some sense, and deleting them from it. 
The term has its origins in the statistical practice of identifying ``outliers" in a sample of i.i.d. random variables, then removing them from a statistic of interest, typically, the sample sum, which can be considered as a stochastic process in discrete time.
More recently, the techniques have been transferred to processes such as extremal processes and L\'evy processes indexed by a
 continuous time parameter, 
where asymptotic properties of the trimmed process have been worked out in a number of interesting cases.
The asymptotic studied may be large time ($t\to\infty$), as in the statistical situation, or, for continuous time processes, small time ($t\dto 0$). The small time case extends our understanding  of local properties of the process and can have direct application as for example in
Maller and Fan \cite{FM2015} and Maller and Schmidli \cite{MS2017}; the large time case has the statistical applications alluded to, such as the robustness of statistics, and insurance modelling, etc., as we discuss later.\\ \indent
This area of research can be regarded as combining studies on properties of extremes of the jumps of a process with those of the process itself; the former, a version of extreme value theory; the latter relating, for example, to domains of attraction of the process. A combination of the two fields enriches both, and the trimming idea is a natural way of approaching this. 
Research in a similar direction has been carried on by  Silvestrov and Teugels \cite{silvestrov2004}; see also their references.
\\ \indent
In the present paper we extend some earlier work of the present authors to consider various ways of trimming the sample path of a stochastic process in the space 
$ \DD[0,1]$ of \cadlag \ functions. The initial set-up is very general. We begin in Section \ref{s1} by establishing continuity properties (in the Skorokhod (strong) $J_1$-topology) of operators which remove extremes.
There are a number of intuitively reasonable ways of defining such operators. Not all of them behave in the same way, and Section \ref{s1} is devoted to teasing out the differences between them. We take a dynamic {\it sample path approach} which brings into focus some interesting and distinctive features not previously apparent. Proofs for Section \ref{s1} are in Section \ref{s4}.

An application of the ideas to the  functional convergence of a L\'evy process in the domain of attraction of a stable law 
 is then given.
 Statements for these are in Section \ref{s2}, and proofs of them are in Section \ref{s3}. 
 Continuity properties of certain extremal operators are closely related to the occurrence or otherwise of tied (equal) values in the large jumps of the limiting process and consequently we need to analyse these too.
 A final Section \ref{S6} develops a motivating application to a reinsurance ruin time problem.


\section{Extremal Operators on Skorokhod Space}\label{s1}
Let $\DD([0,1], \R) = : \DD$ be the space of all  \cadlag \, functions: $[0,1] \to \R$  endowed with the Skorokhod (strong) $J_1$-topology.
Denote the sup norm by $||\cdot||$, so that $||x||=\sup_{0\le \tau\le 1}|x(\tau)|$, where
for each $x\in \DD$, $x(\tau)$ is the value of $x$ at time $\tau\in[0,1]$. 
Convergence in the $J_1$-topology is characterised as follows.
Let $\Lambda$  be the set of all continuous and strictly increasing functions $\lambda : [0,1] \to [0,1]$  with $\lambda(0) = 0$ and $\lambda (1) = 1$. 
Denote by $I: [0,1] \to [0,1]$ the identity map.
Let $\alpha_n \in \DD$. Then $\alpha_n \toj \alpha$ in $\DD$ if there exists a sequence $(\lambda_n) \in \Lambda $ such that 
\ben\label{J1cvg}
||\lambda_n - I|| \vee || \alpha_n \circ \lambda_n - \alpha || \to 0  \quad \text{as}\quad n \to \infty.
\een
Intuitively,  $J_1$-convergence requires ``matching jumps'' at ``matching points'' after a deformation of time.
We refer to Chapter VI in Jacod and Shiryaev \cite{JacodShiryaev03} and Section 12 in Billingsley \cite{bill1999} for more information on the Skorokhod space. For other topologies on $\DD$ we refer to Skorokhod \cite{sko56} and Whitt \cite{whitt2002}.

We proceed by setting out some basic methods of trimming extremes.

\subsection{Global (Pointwise) Trimmers (``Trim As You Go'')}\label{basic}
 Let $x =(x(\tau))_{0\le \tau\le 1} \in \DD$ with jump process $\Delta=(\Delta  x(\tau):=x(\tau)-x(\tau-))_{\tau>0}$. 
Set $\Delta x(0)\equiv 0$.  Define the following extremal operators mapping
 $ \DD$ into $\DD$:

(i)\ $\SSS(x) (\cdot)  = \sup_{0\le s\le \cdot} x(s)$;

(ii)\ $\wt \SSS (x) (\cdot)  = \sup_{0\le s\le \cdot} |x(s)|$;

(iii)\  $ \SSS_{\pm\Delta} (x) (\cdot)  = \SSS_{\Delta} (\pm x) (\cdot) :=  \sup_{0\le s\le \cdot} \Delta (\pm x)(s) \wedge 0 $; and

(iv)  $ \wt \SSS_{\Delta} (x) (\cdot)  = \sup_{0\le s\le \cdot} |\Delta x(s)|$.

\noindent 
Here (i)  $\SSS(x)$ is the running supremum process of $x$,
and  (ii)  $\wt \SSS(x)$  is the running supremum process for $|x|$.
In (iii),  $ \SSS_{+\Delta} (x) (\tau)$ represents the magnitude of largest (positive) jump of $x$ up till time $\tau$, and  $ \SSS_{-\Delta} (x) (\tau)$
represents the largest magnitude of the negative jumps of $x$ up till time $\tau$; while in (iv),   $ \wt \SSS_{\Delta} (x) (\tau)$
represents the magnitude of largest jump in modulus of $x$ up till time $\tau$.

With these operators we can define what we call {\it global, or pointwise, trimming operators}.
 Let $\NN:=\{1,2,\ldots\}$, $\NN_0:=\NN\cup\{0\}$.
Take $r=2,3,\ldots$ and define iteratively

(v)\  the $r^{th}$ extremal positive (negative) trimming operators
 \[
 \TTT^{(1, \pm)}_{rim} (x)  = x \mp \SSS_{\pm\Delta} (x)  \quad \text{and} \quad 
 \TTT_{rim}^{(r,\pm)}(x)   = \TTT^{(1,\pm)}_{rim} \circ \TTT_{rim}^{(r-1,\pm)} (x);
\]

(vi)\ the $r^{th}$ extremal positive (negative) jump operators
\[
\SSS^{(1, \pm)}_{\Delta} (x)  = \SSS_{\pm\Delta} (x)  
\quad \text{and} \quad 
 \SSS^{(r, \pm)}_{\Delta} (x)  = \SSS^{(1, \pm)}_{\Delta} \circ \TTT_{rim}^{(r-1, \pm)} (x);
\]

(vii)\  the $r,s$ trimming  operators (for $s\in\NN$)
 \[
 \TTT_{rim}^{(r,s)}(x)   = \TTT^{(r, +)}_{rim}  \circ \TTT^{(s,-)}_{rim}(x)  = \TTT^{(s, -)}_{rim} \circ \TTT^{(r, +)}_{rim} (x);
\]

(viii)\ the $r^{th}$ extremal modulus trimming  operators
 \[
\wt \TTT_{rim} (x) = \wt \TTT^{(1)}_{rim} (x) = x - \wt \SSS_{\Delta} (x)  \quad \text{and} \quad 
  \wt \TTT_{rim}^{(r)}(x)   = \wt \TTT_{rim} \circ \wt \TTT_{rim}^{(r-1)} (x); 
\]

(ix)\ and the $r^{th}$ extremal modulus jump operators
 \[
\wt \SSS^{(1)}_{\Delta} (x)  = \wt \SSS_{\Delta} (x)  \quad \text{and} \quad 
  \wt \SSS^{(r)}_{\Delta} (x)  = \wt  \SSS_{\Delta} \circ \wt \TTT^{(r-1)}_{rim} (x).
\]

\noindent
Here  (v) $\TTT_{rim}^{(r,+)}(x)(\tau)$ is $x$ with the $r$ largest jumps of $x$ up till time $\tau$ subtracted, and
 $\TTT_{rim}^{(r,-)}(x)(\tau)$ is similar with the $r$ negative jumps of largest magnitude subtracted.
 In (vii),  $\TTT_{rim}^{(r,s)}(x)$ has the $r$ positive and $s$ negative jumps of largest magnitudes subtracted, while (viii), 
 $\wt \TTT_{rim}^{(r)}(x)(\tau)$ has  the $r$ largest jumps in modulus of $x$ up till time $\tau$ removed from $x$.  In (vi) and (ix), 
 $ \SSS^{(r, \pm)}_{\Delta} (x)$ and $  \wt \SSS^{(r)}_{\Delta} (x)$
 are the $r^{th}$ largest values in magnitude for positive (negative), or in modulus, jumps of the corresponding processes.
  
 We call the  operators in (v), (vii), (viii), ``Trim As You Go'' operators because at each point in time, the designated number of  largest positive (negative) jumps up to that point are removed from the process.
See Figure \ref{trimasyougo} in Section \ref{location} for an illustration with a schematic insurance risk process. 

To analyse the convergence of these operators in $\DD$, we need the following considerations.
We say that an operator  $\Psi:\DD\to \DD$ is {\it $\|\cdot\|$-continuous} at $x\in\DD$ if $\lim_{n\to\infty}\|x_n-x\|= 0$ implies $\lim_{n\to\infty}\|\Psi(x_n)-\Psi(x)\|= 0$.
We say that  $\Psi$ is {\it $J_1$-continuous} at $x\in\DD$ if $x_n\toj x$ implies $\Psi(x_n)\toj \Psi(x)$, as $n\to\infty$. 
In general,  $\Psi$ being $\|\cdot\|$-continuous at $x$
does not imply that  $\Psi$ is $J_1$-continuous at $x$. 
However, this is true if in addition $\Psi$ is {\it $\Lambda$-compatible}, by which we mean $\Psi(x)\circ \lambda=\Psi(x\circ\lambda)$ for all $x\in\DD$, $\lambda\in\Lambda$. 
The operator $\Psi$ is called {\it jointly $J_1$-continuous at $x$} if for any sequence $x_n$ converging
to $x$ in the
$J_1$-topology, there exists
  $(\lambda_n)$ in $\Lambda$, such that, simultaneously
as $n \to \infty$, $\|\lambda_n - I\| \to 0$, $\|x_n \circ \lambda_n - x\| \to 0$ and $\|\Psi(x_n) \circ \lambda_n - \Psi(x)\| \to 0$.
The following simple proposition summarises. 

\begin{prop}\label{propskounif}
Let $\Psi:\DD\to \DD$ be $\Lambda$-compatible and take $x\in\DD$. 
Suppose  $\Psi$ is $\|\cdot\|$-continuous at $x$.
  Then $\Psi$ is jointly $J_1$-continuous at $x$.
\end{prop}

\noindent   
{\bf Proof of Proposition~\ref{propskounif}:}\ 
Assume $\Psi$ is $\Lambda$-compatible and  $\|\cdot\|$-continuous at $x\in\DD$.
Since  $\Psi(x_n)\circ \lambda_n=\Psi(x_n\circ \lambda_n)$, 
$x_n\toj x$, i.e.,  $\lim_{n\to\infty}\|x_n\circ\lambda_n-x\|=0$,
together with   $\Psi$ being $\|\cdot\|$-continuous at $x$, 
 implies $\lim_{n\to\infty}\|\Psi(x_n)\circ\lambda_n-\Psi(x)\|=0$, 
 i.e., $\Psi(x_n)\toj \Psi(x)$,  proving the proposition. 
\halmos


\begin{prop}\label{operCon}
Each of the operators defined in {\rm (i)--(ix)} is
\begin{enumerate}[\rm (a)]
\item  $\|\cdot\|$-Lipschitz, hence continuous in $\|\cdot \|$ norm;
\item  $\Lambda$-compatible; and, consequently,  by  Proposition \ref{propskounif}, jointly $J_1$-continuous.
\end{enumerate}
\end{prop}

\noindent   {\bf Proof of Proposition \ref{operCon}:}\ 
(a)\ For example, we prove (iii).
When $x,y\in\DD$, 
\bea\label{2in}
&&\|\SSS_\Delta(x)-\SSS_\Delta(y)\|
=
\sup_{0<t\le 1}\|\sup_{0<s\le t}(x(s)-x(s-))-\sup_{0<s\le t}(y(s)-y(s-))\|\cr
&&\cr
&\le&
\sup_{0<s\le 1}\|(x(s)-x(s-))-(y(s)-y(s-))\|
\le
2\|x-y\|,
\eea
using  the triangle inequality
$ \big| ||\alpha|| - ||\beta||      \big| \le ||\alpha- \beta|| \le ||\alpha|| + ||\beta||$.

\noindent
{\rm (b)}\
We prove this for (iv),  for example.
Let $x \in \DD$, $\lambda \in \Lambda$, and $t \in [0,1]$.
Then 
\[ 
\wt \SSS_\Delta(x \circ \lambda) (t) = \sup_{0\le s\le t}| \Delta (x\circ \lambda)(s)| = \sup_{0\le s\le t} |\Delta x(\lambda(s))| = \sup_{0\le s\le \lambda(t)} |\Delta x(s)| = \wt \SSS_\Delta(x) (\lambda(t)).
\]
\vskip-0.3cm
\halmos

\begin{remark}
{\rm 
We refer to Section IV.2 in Jacod and Shiryaev \cite{JacodShiryaev03} for other continuity properties of common mappings in the Skorokhod topology.
}
\end{remark}

\subsection{Signed Modulus Trimmers}\label{SMT}
In Section \ref{s2} we will consider a L\'evy process $X=(X_t)_{t\ge 0}$ which is to be trimmed.
Before this, in the present subsection, we want to draw attention to an issue that arises with modulus trimming when considered pathwise.
There may be one or more jumps equal in magnitude to the largest of $|\Delta X_s|=|X_s-X_{s-}|$, for $0\le s\le t$. We refer to these as ``tied" values (for the modulus, with a similar concept for the positive and negative jumps).

Buchmann, Fan and Maller \cite{bfm2016} (hereafter, ``BFM") 
define  a ``modulus trimmed L\'evy process" as follows. 
Denote the largest modulus jump of $X$ up to time $t$, i.e., the jump corresponding to the largest of $|\Delta X_s|$, $0\le s\le t$,
  by $\wt{\Delta X}_t^{(1)}$. 
 When there is no tie for $\sup_{0\le s\le t}|\Delta X_s|$,
 the sign of  $\wt {\Delta  X}_t^{(1)}$ is uniquely determined.
 When there is a tie, the procedure in BFM is to nominate a jump 
chosen at random  among the almost surely (a.s.) finite number of tied values according to a discrete uniform distribution on the collection of ties. While appropriate in the context of BFM, 
 this definition is problematic when we consider the sample path of the process on $[0,1]$.
 To see why, take a simple example. 
 Suppose for some $\omega$ the largest modulus jump up till time $t$ is tied at values $0< s_1 < s_2 <t$ with opposite signs, say:
 \[  \Delta X_{s_1}(\omega) = |\wt{\Delta X}^{(1)}_t (\omega)| \quad \text{and} \quad  \Delta X_{s_2}(\omega) = -|\wt{\Delta X}^{(1)}_t (\omega)|,
 \]
while $ |\wt{\Delta X}^{(1)}_{s}(\omega)| = |\wt{\Delta X}^{(1)}_t (\omega)|$ for all $s \in [s_2, t]$.
For each $s \in [s_2,  t]$, if we were to choose from $\{\Delta X_{s_1}, \Delta X_{s_2}\}$ with equal probability to be trimmed from $X_s(\omega)$, the sample path of the resulting trimmed process would  not be in $\DD$.

Thus, we need to design a way to define signed modulus trimming on the sample path of $X$ so as to stay within $\DD$.  
 One way to do this is as follows (we now revert to the general setup).
For $x\in\DD$, define the {\it last modulus record time process} on $[0,1]$ as
\be\label{Ltdef}
\wt L_\tau(x):=\sup \{ s \in[0,\tau]: |\Delta x(s)| = \wt \SSS_\Delta(x)(\tau)\}, \ {\rm for\ each}\  \tau \in [0,1].
\ee
Then  the {\it signed largest modulus jump} up till time $\tau\in [0,1]$ is
$ \Delta x(\wt L_\tau(x))$,
and the {\it signed largest trimmer} can be defined as
$\TTTT_{rim}(x):= x- \Delta x(\wt L_\tau(x))$.
More generally, interpret $\TTTT^{(0)}_{rim}(x) = x$, let
$\TTTT^{(1)}_{rim}(x) := \TTTT_{rim}(x)$, and, for $r=2,3,\ldots$, set
\ben\label{def_mod}
\TTTT_{rim}^{(r)}(x):= \TTTT_{rim}( \TTTT_{rim}^{(r-1)}(x)).
\een

Now $\TTTT^{(1)}_{rim}$ is not in general globally 
$J_1$-continuous, that is, it is not  $J_1$-continuous at all $x\in \DD$.
Take for example $x = {\bf 1}_{[1/3, 2/3)}$ and $x_n = x + \frac 1 n {\bf 1}_{[1/3, 1]}$. Then $x_n \toj x$, but 
\[  \TTTT_{rim}^{(1)} (x) = {\bf 1}_{[2/3,1]} \quad \text{and} \quad  \TTTT_{rim}^{(1)} (x_n) =  -{\bf 1}_{[2/3,1]}.
\] 
However, $\TTTT_{rim}$  is continuous when there is no change of signs of ties in the limit. This is shown in Theorem \ref{signTies}, which uses   the following notation. 
For each $\tau \in [0,1]$, collect the times of occurrence of the largest values, and the times of occurrence of values having largest modulus, into sets 
 $\myAA^\pm_\tau(x)$ and $\widetilde \myAA_\tau(x)$, thus:
\be\label{Apm}
\myAA^\pm_\tau(x):=\{0\!<\!s\!\le\! \tau:\Delta x(s)=\SSS_{\pm \Delta}(x)(\tau)\}
\ee
 and
\be\label{Atil}
\widetilde \myAA_\tau(x):=\{0\!<\! s\!\le\! \tau:|\Delta x(s)|\!=\!\widetilde \SSS_\Delta(x)(\tau)\}.
\ee
We use the convention that when $x$ is continuous on $[0,\tau]$, then $\myAA_\tau^\pm(x) = \wt \myAA_\tau (x)=\emptyset$.
Recall that a  \cadlag \, function has only finitely many jumps with magnitude bounded away from 0, so $\myAA_\tau^\pm(x)$ and $\widetilde\myAA_\tau(x)$ are 
finite sets (we include in this the possibility that one or other of them may be empty) for functions $x\in\DD$. Collect the sign changing largest modulus jumps contained  in $ \wt \myAA_\tau(x) = \{s_1, \ldots s_{\# \wt \myAA_\tau(x)}\}$ into the set 
\[ 
\BB_\tau(x) : = \{ s_k \in \wt \myAA_\tau(x) : \Delta x(s_{k}) = - \Delta x(s_{k-1}), \text{ where } k = 2, \ldots, \# \wt \myAA_\tau(x)  \}.
\]
Note that $ \# \wt \myAA_\tau(x) = 1$ implies $\# \BB_\tau(x) = 0$. Conversely, $\# \BB_\tau(x) = 0$ implies $ \wt \myAA_\tau(x) = \myAA^+_\tau(x)$ or 
$ \wt \myAA_\tau(x) = \myAA^-_\tau(x)$.

In the next theorem, we show that when there is no sign change among ties of large modulus jumps in $x$ and its trimmed versions $\TTTT_{rim}^{(j)}(x)$ for all $0\le j \le r-1$,  $\TTTT_{rim}^{(r)}$ is jointly $J_1$-continuous at $x$. 
The next theorem is proved in Section \ref{s4}.
\begin{theorem}\label{signTies} 
$\TTTT_{rim}$ is jointly $J_1$-continuous at $x$ if $\sup_{\tau\in [0,1]}\#\BB_\tau(x) = 0$. Consequently,  $ \TTTT^{(r)}_{rim}$ is jointly $J_1$-continuous at $x$ if $\sup_{\tau\in [0,1]}\#\BB_\tau(\TTTT^{(j)}_{rim}(x)) = 0$ for all $j=0,\ldots, r-1$,
when $r\in\NN$. 
\end{theorem}

\subsection{Record Times Trimmers (``Lookback Trimming'')}\label{location}
\begin{figure}[h!]
\centering
\includegraphics[width=.43\linewidth]{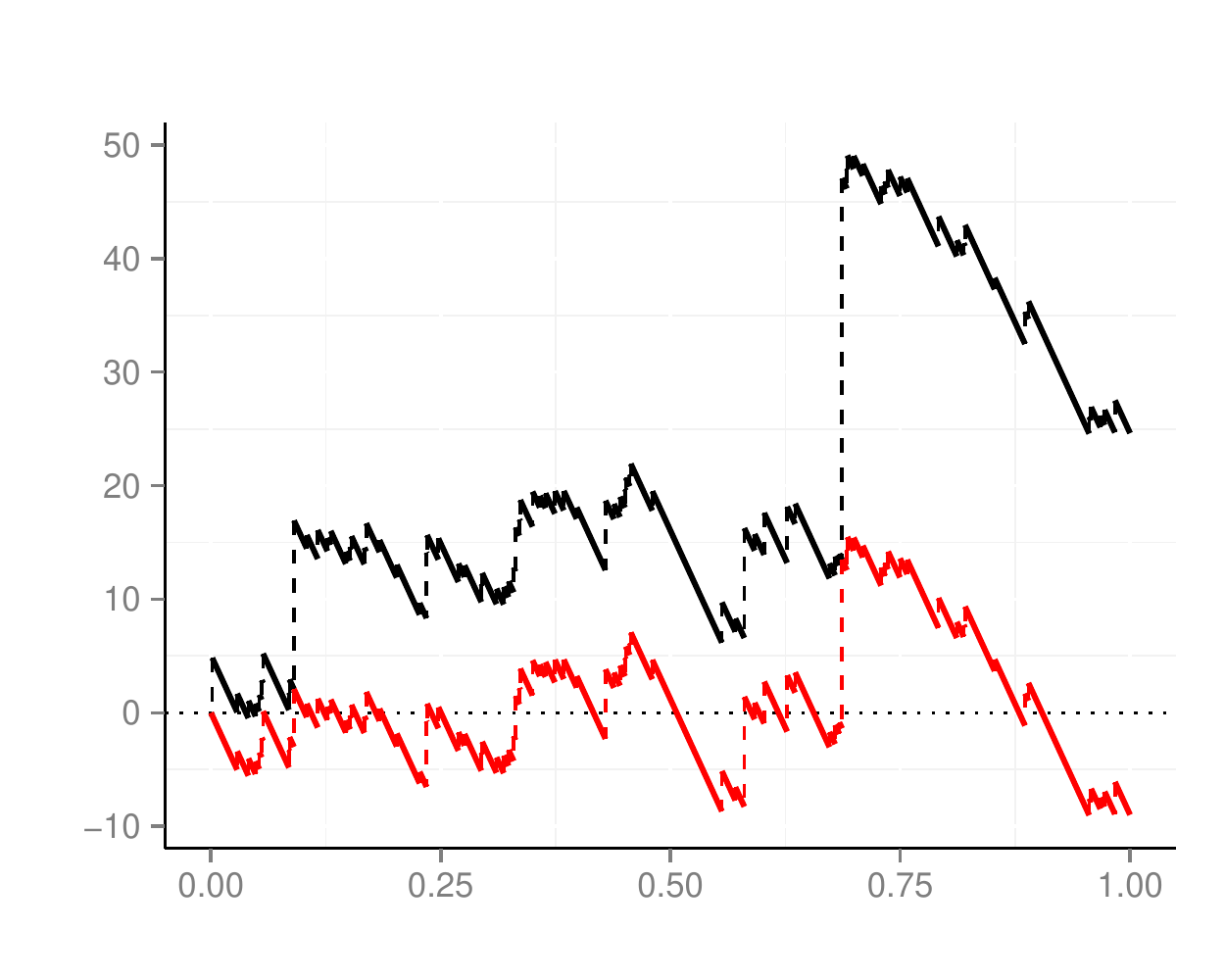}
\includegraphics[width=.43\linewidth]{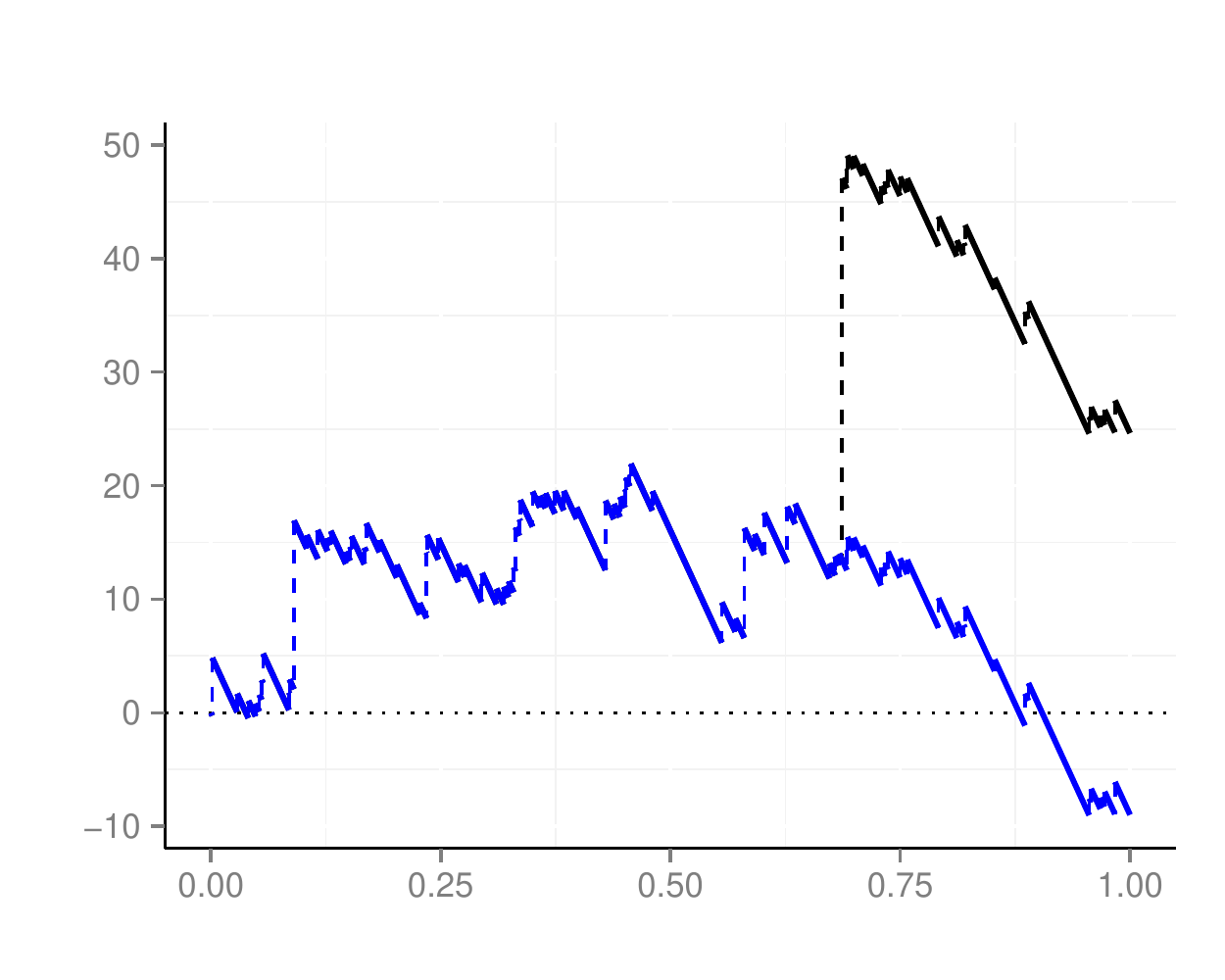}
\caption{\footnotesize A realisation of the compound Poisson risk insurance model $X_t = ct + \sum_{i=1}^{N_t} \xi_i$, $t \in [0,1]$, is represented in black. The $\xi_i$ are i.i.d. with Pareto(1,2) distribution, $c= -110$ and $N_t$ is Poisson with $\EE(N_1) = 100$. Left: the sample path of $\TTT_{rim}^{(1,+)}((X_t)_{t \in [0,1]})$ (``trim as you go'')  is represented in red. Right: the sample path of $\RRR_{trim}((X_t)_{t \in [0,1]})$ (``lookback trimming'') is represented in blue. This path coincides with the original process before $ R_1((X_t)_{t \in [0,1]})$ and with the $\TTT_{rim}^{(1,+)}((X_t)_{t \in [0,1]})$ path after.}
\label{trimasyougo}  
\end{figure}

On the function space $\DD$, we can extend the idea of trimming by including a random location where trimming starts
and hence define a second kind of trimming.
For $x\in\DD$ define {\it the first (positive) record time} in $[0,1]$ by
\be\label{Rdef}
R_\tau(x):=\inf\{s\in[0,\tau]:\Delta x(s)=\SSS_{+\Delta}(x)(\tau)\}
\ee
for $0<\tau\le 1$, and similarly we could define the first (negative) record time. Likewise, 
\be\label{Rtdef}
\wt R_\tau(x):=\inf\{s \in[0,1]:|\Delta x(s)|=\widetilde \SSS_\Delta(x)(\tau)\}
\ee
gives the \emph{first modulus record time}. The corresponding {\it record time trimmers} are
\be\label{defRtr}
\RRR_{trim}(x):= x-\Delta x(R_1(x)) {\bf 1}_{[R_1 (x),1]},\quad
\widetilde \RRR_{trim}(x):= x-\Delta x( \wt R_1(x)) {\bf 1}_{[\wt R_1(x),1]}.
\ee
Expanding, $\RRR_{trim}(x)(\tau)$ can be written for $\tau\in[0,1]$ as
\ben
\RRR_{trim}(x)(\tau)=
\begin{cases} x(\tau)-\sup_{0<s\le \tau}\Delta x(s), &\mbox{if } R_1(x)\le \tau;  \\
x(\tau), & \mbox{otherwise. } 
\end{cases} 
\een
Thus, $x$ is trimmed at time $\tau$ if the record occurs before $\tau$, otherwise not. 
Figure \ref{trimasyougo} gives an illustration of the two trimming types for a compound Poisson risk process as used 
in insurance risk modelling (cf. Section \ref{S6}).
 
Set $\RRR_{trim}^{(0)}(x) = x$. For $r\in\NN$,  define the \emph{$r^{th}$ order record times trimmers} as 
\be\label{defRt}
\RRR_{trim}^{(r)}(x):= \RRR_{trim}(\RRR_{trim}^{(r-1)}(x)) =\RRR_{trim}^{(r-1)}(x) -
\Delta \RRR_{trim}^{(r-1)}\big( R_1 \circ \RRR_{trim}^{(r-1)}(x)\big)
{\bf 1}_{[R_1 \circ \RRR_{trim}^{(r-1)}(x),1]},
\ee
and, for $\wt \RRR_{trim}^{(r)}(x)$, replace each $R_1$, $\RRR_{trim}$ in \eqref{defRt} by $\wt R_1$ and $\wt \RRR_{trim}$. 


While the record trimming functionals are $\Lambda$-compatible, they are not however of the type  described in Section \ref{basic}. In fact they are not  globally  $J_1$-continuous.

\begin{example} \label{egrtrim}{\bf [$\RRR_{trim}$ is  not globally norm or $J_1$-continuous.]}

(i)\ $\RRR_{trim}$ is  not $\|\cdot\|$-continuous.
To see this,  let $x:={\bf 1}_{[1/3,1]} +{\bf 1}_{[2/3,1]}$. 
For $n\in\NN$ set $x_n:=x+(1/n) {\bf 1}_{[2/3,1]}$.
Observe that $\lim_{n\to\infty}\|x_n-x\|\le \lim_{n\to\infty}1/n=0$. In particular, $x_n\sko x$ as $n\to\infty$. 
However, 
\[R_1(x)=1/3,\; \RRR_{trim}(x)={\bf 1}_{[2/3,1]},\
R_1(x_n)=2/3,\; \RRR_{trim}(x_n)={\bf 1}_{[1/3,1]},\quad n\in\NN,\]
and
\[\|\RRR_{trim}(x_n)-\RRR_{trim}(x)\|=1,\ n\in\NN.\]
Thus, $\RRR_{trim}$ is $\Lambda$-compatible, but not $\|\cdot\|$-continuous in $x$.

(ii)\
 $\RRR_{trim}$ also is not $J_1$-continuous. To see this, take  $(\lambda_n)\subseteq\Lambda$ with $\lim_{n\to\infty}\|\lambda_n-I\|=0$. Then 
$\lim_{n\to\infty}\lambda_n^{-1}(1/3)=1/3$, hence 
for $n\in\NN$, once $\lambda_n^{-1}(1/3)<2/3$, 
\begin{eqnarray*}
\|\RRR_{trim}(x)-\RRR_{trim}(x_n\circ \lambda_n)\|&=&\|\RRR_{trim}(x)-\RRR_{trim}(x_n)\circ \lambda_n\|\\
&\ge& |\RRR_{trim}(x)(\lambda_n^{-1}(1/3))-\RRR_{trim}(x_n)(1/3)|\\
&=& {\bf 1}_{[1/3, 2/3]}(1/3)= 1. \end{eqnarray*}
Consequently, we have $x_n\sko x$, but not $\RRR_{trim}(x_n)\sko \RRR_{trim}(x)$. \halmos
\end{example}

Recall the definitions of $\myAA^\pm_\tau$ and $\wt \myAA_\tau$ in \eqref{Apm} and \eqref{Atil}. Our main result of this section is that the record time trimmer $\RRR_{trim}$ is jointly $J_1$-continuous at $x$ if and only if $x$ does not admit ties.
The next theorem is proved in Section \ref{s4}.

\begin{theorem}\label{theotrimcont} Let $x\in\DD$ and $r\in\NN$.\\
{\rm (i)} If $\#\myAA^+_1(x)\le 1$ then $\RRR_{trim}$ is jointly $J_1$-continuous at $x$. Consequently, if \\
\mbox{$\#\myAA^+_1(\RRR_{trim}^{(j)} (x))\le 1$} for all $j = 0, \ldots, r-1$, then $\RRR_{trim}^{(r)}$ is jointly $J_1$-continuous at $x$.
\\[1mm]
{\rm (ii)} If  $\RRR_{trim}$ is $J_1$-continuous at $x$ then $\#\myAA^+_1(x)\le 1$.\\
The same holds true with $\RRR_{trim}$ replaced by $\wt \RRR_{trim}$ and $\#\myAA^+_1(x)$ replaced by $\#\wt \myAA_1(x)$.
\end{theorem}

\section{Functional Laws for L\'evy Processes}\label{s2}
A number of interesting processes can be derived by applying the operators in Section \ref{s1} to L\'evy processes. 
In the present section $X=(X_t)_{t>0}$, $X_0=X_{0-}=0$,  will be a real valued \cadlag \, L\'evy process with canonical triplet $(\gamma, \sigma^2=0, \Pi_X)$. 
%
%
%
The positive, negative and two-sided tails of the L\'evy measure $\Pi_X$ are, for $x>0$, 
\ben\label{pidef}
\pibar_X^+(x):= \Pi_X\{(x,\infty)\},\ \pibar_X^-(x):= \Pi_X\{(-\infty,-x)\},\
{\rm and}\  \pibar_X(x):=\pibar_X^+(x)+\pibar_X^-(x).
\een

The jump process of $X$ is $(\Delta X_t = X_t-X_{t-})_{t\ge 0}$, the positive jumps are $\Delta X_t^+ = \Delta X_t \vee 0$, and the magnitudes of the negative jumps are $\Delta X_t^- = (-\Delta X_t)\vee 0$.
The processes $(\Delta X_t^+)_{t\ge 0}$ and $(\Delta X_t^-)_{t \ge 0}$, when present,  are non-negative independent processes. 
For any integers $r,s  > 0$, let $\Delta X_t^{(r)}$ be the $r^{th}$ largest positive jump, and let $\Delta X_t^{(s),-}$ be the 
$s^{th}$ largest jump in $\{\Delta X_s^{-}, 0<s\le t\}$, 
i.e., the negative of the $s^{th}$ smallest  jump. These kinds of ordered jumps are carefully defined in BFM, allowing for the possibility of tied values. (Recall the discussion in Subsection \ref{SMT}.)
We can similarly define  $\Delta X_{t-}^{(r)}$ and $\Delta X_{t-}^{(s),-}$ for the ordered jumps in $\{\Delta X_s, 0<s<t\}$.

Throughout, for small time convergence ($t\dto 0$) we assume $\pibar_X(0+) = \infty $ when dealing with modulus trimming and $\pibar_X^+(0+) = \infty$ or $\pibar_X^-(0+) = \infty$ (or both when appropriate) when dealing with one-sided trimming. In particular, these ensure
there are infinitely many jumps $\Delta X_t$, or $\Delta X_t^\pm$,  a.s., in any bounded interval of time. 

As  demonstrated in Section \ref{s1}, the largest modulus trimming as defined in BFM is not natural for the functional setting, so here we adopt a modified definition.
Write $\wt{\Delta X}_t^{(r)}$ to denote the $r^{th}$ largest jump in modulus up to time $t$, taking the sign of the {\it latest} $r^{th}$ largest modulus jump. Then define the trimmed L\'evy processes  
\begin{equation}\label{TL}
{}^{(r, s)} X_t:= X_t- \sum_{i=1}^r {\Delta X}_t^{(i)} + \sum_{j=1}^s \Delta X_t^{(j),-}
 \quad {\rm and} \quad
{}^{(r)}\wt X_t:= X_t- \sum_{i=1}^r \wt{\Delta X}_t^{(i)},
\end{equation} which we call the {\it asymmetrically trimmed} and 
{\it modulus trimmed} processes, respectively. 
With the convention $\sum_1^0\equiv 0$, taking $r = 0$ or $s = 0$ in asymmetrical trimming gives {\it one-sided trimmed processes}
$ {}^{(r)} X_t:={}^{(r,0)} X_t$
and ${}^{(s,-)} X_t:= {}^{(0, s)} X_t$.

In this section we apply a functional law for L\'evy processes  attracted to a non-normal stable law to get two theorems for trimmed L\'evy processes.
$X_t$ is said to be in a {\it non-normal domain of attraction} 
at small (large) times if there exist non-stochastic functions $a_t \in \R$ and $b_t > 0$ such that 
\be\label{dom}
 \frac{X_t - a_t}{b_t} \todr Y \quad \text{as } t \dto 0\, (t\to\infty),
\ee 
where $Y$ is an a.s. finite, non-degenerate, non-normal\footnote{When $Y$ is $N(0,1)$, a standard normal random variable, the large jumps are asymptotically negligible with respect to $b_t$ and \eqref{ffunLook} and \eqref{ffunLaw} remain true with $^{(r,s)}\YYY$ 
and  $^{(r)}\wt\YYY$ a standard Brownian motion; see Fan \cite{fan2015JOTP}.}  random variable. 
Then \eqref{dom} implies that the two-sided tail $\pibar_X$ of $X$ is regularly varying (at 0 or $\infty$, as appropriate),  with index 
$\alpha\in(0,2)$.  The limit rv $Y$ has the distribution of  $Y_1$,  where $(Y_\tau)_{0\le \tau \le 1} \equiv \YYY$ is a stable($\alpha$) L\'evy process.
The canonical triplet for $\YYY$ will be taken as $(0, 0, \Pi_Y)$, where $\Pi_Y$ has tail function $\pibar_Y(x) = cx^{-\alpha}$, $x>0$, for some $c >  0$. 

In the small time case, conditions on the L\'evy measure for \eqref{dom} to hold can be deduced from Theorem 2.2 of Maller and Mason \cite{MM2010}, whose result can also be used to show that \eqref{dom} can be extended to convergence in $\DD$;
that is, 
\be\label{funcFC1}
\III_t =  \{\III_t(\tau)\}_{0\le \tau \le 1} :=
\left( \frac{X_{\tau t} - \tau a_t}{b_t} \right)_{0\le \tau \le 1}
\to (Y_\tau)_{0\le \tau \le 1} = \YYY,
\ee
weakly as $t \dto 0$ with respect to the $J_1$-topology. Large time ($t \to \infty$) convergence in \eqref{funcFC1} also follows from \eqref{dom} as is well known. 

Assuming the  convergence in \eqref{funcFC1}, we can prove a variety of interesting functional limit theorems for $X$ by applying the operators in Section \ref{s1}.  We list some examples in  Theorem \ref{maxtrim}  and prove them in Section \ref{s3}.

  Theorem \ref{maxtrim}  considers (i) lookback trimming, (ii) two-sided (or one-sided, with $r$ or $s$ taken as 0) trimming, and (iii) signed modulus trimming defined as in Subsection \ref{SMT}.
To specify the lookback trimming in this situation, recall the 
definition of the record time trimming functionals in 
 \eqref{defRtr} and \eqref{defRt}.
 Using them, we define, for $X$,
 {\it lookback trimmed paths of order $r$}, based on {\it positive} jumps, being processes on $\tau\in[0,1]$, indexed by $t >0$, as
\[
(^{(1)}X^R_{t\tau})_{\tau\in [0,1]} = \RRR_{trim}\big((X_{t\tau})_{\tau\in [0,1]}\big),
\]
 and, for $r = 2,3, \ldots,$
\[
 (^{(r)}X^R_{t\tau})_{\tau \in [0,1]} = \RRR^{(r)}_{trim}\big((X_{t\tau})_{\tau \in [0,1]}\big) =\RRR_{trim}\big(\RRR^{(r-1)}_{trim}\big((X_{t\tau})_{\tau\in [0,1]}\big)\big).
\]

\begin{theorem}\label{maxtrim}
Assume $(X_t)_{t\ge 0}$ is in the domain of attraction of a stable law at $0$ (or $\infty$) with nonstochastic centering and norming  functions $a_t\in\R$, $b_t>0$, so that \eqref{dom} and \eqref{funcFC1} hold.
In the following, convergences are with respect to the $J_1$-topology
in $\DD$. 

{\rm (i)}\ Suppose $\pibar_X^+(0+) = \infty$ and  $r \in \NN$. Then, for the same $a_t$ and $b_t$, 
\be\label{ffunLook}
\left(\frac{{}^{(r)}X^R_{\tau t}-\tau a_t}{b_t}\right)_{ 0 \le \tau \le 1} \todr ({}^{(r)}Y_\tau^R)_{0\le \tau \le 1} :=  {}^{(r)}\YYY^R \ \text{in} \ \DD, \ \text{as} \  t \dto 0 \ (t \to \infty).
\ee 

{\rm (ii)}\ 
Assume $\pibar_X^+(0+)=\pibar_X^-(0+) = \infty$ and $r,s \in \NN_0$.  Then, for the same $a_t$ and $b_t$, 
\be\label{ffunLaw}
\left(\frac{{}^{(r,s)}X_{\tau t}-\tau a_t}{b_t}\right)_{ 0 \le \tau \le 1} \todr \big({}^{(r,s)}Y_\tau\big)_{0\le \tau \le 1} := {}^{(r,s)}\YYY \ \text{in} \ \DD, \  \text{as} \  t \dto 0 \ (t \to \infty).
\ee

{\rm (iii)}\ 
Suppose only $\pibar_X(0+) = \infty$ and  $r\in\NN$.
 Then, for the same $a_t$ and $b_t$, 
 \be\label{ffunLaw2}
\left(\frac{{}^{(r)}\wt X_{\tau t}-\tau a_t}{b_t}\right)_{ 0 \le \tau \le 1} \todr \big({}^{(r)}\wt Y_\tau\big)_{0\le \tau \le 1} :={}^{(r)}\wt \YYY  \ \text{in} \ \DD, \  \text{as} \  t \dto 0  \ (t \to \infty).
\ee 
\end{theorem}

\begin{remark}
{\rm 
(i)\
 BFM and also Maller and Mason \cite{MM2010} include convergence of the~quadrat\-ic variation of $X$ in their expositions. Using these as basic convergences (i.e., together with \eqref{dom} and \eqref{funcFC1})  would lead to functional convergences of the jointly trimmed process together with its trimmed quadratic variation process, and we could then consider self-normalised versions. But we omit the details of these.\\[1mm] 
(ii)\ Fan \cite{fan2015SPA} proved the converses in (ii) and (iii) of  Theorem \ref{maxtrim} for $t \dto 0$,
i.e., if the convergence in \eqref{ffunLaw} or  \eqref{ffunLaw2} holds for a fixed $\tau>0$, then
$X$ is in the domain of attraction of a stable law with index $\alpha\in(0,2)$ at small times. 
}
\end{remark}

We conclude this section by mentioning that the same methods can be used to get functional convergence for jumps of an extremal process together with trimmed versions. Again, we omit further details.

\section{Proofs for Section \ref{s1}}\label{s4}
For the  proof of Theorem \ref{signTies} we need a preliminary lemma. 
Let $\CCCC(r_n)$ denote the
set of accumulation points of a sequence $(r_n)\subseteq\RR$ as $n\to\infty$.
Recall that $\wt L_\tau$ and  $\wt R_\tau$,  the last and first modulus record time processes, are defined in \eqref{Ltdef} and \eqref{Rtdef}; 
$R_\tau$, the first positive record time process,  is defined in \eqref{Rdef}. 

\begin{lemma} \label{lemaccpointsT}
Take  $x\in\DD$ and suppose  $(x_n)\subseteq\DD$ with $x_n\sko x$. Then for each $\tau \in [0,1]$,

{\rm (i)}\
if $\wt \myAA_\tau(x)\neq \emptyset$, then $\CCCC(\wt R_\tau(x_n))\subseteq \wt \myAA_\tau(x)$ and $\CCCC(\wt L_\tau(x_n))\subseteq \wt \myAA_\tau(x)$;

{\rm (ii)}\
if $\myAA_\tau^+(x)\neq \emptyset$, then $\CCCC(R_\tau(x_n))\subseteq \myAA^+_\tau(x)$.

{\rm (iii)}\ If $\|x_n-x\| \to 0 $ and  $\#\wt \myAA_\tau(x)=1$, then for all sufficiently large $n$, $\wt R_\tau(x_n) = \wt L_\tau(x_n) = \wt R_\tau(x) = \wt L_\tau(x)$ for each $\tau \in [0,1]$.  
\end{lemma}

 \noindent {\bf Proof of  Lemma \ref{lemaccpointsT}:}\ 
(ii)\
We consider the case $ \myAA^+$ only;
$\wt \myAA$ can be argued similarly.
Take  $x\in\DD$ and let $(x_n)\subseteq\DD$ with $x_n\sko x$. Then there are $\lambda_n\in\Lambda$ such that
$\|\lambda_n-I\|\vee \|y_n-x\|\to 0$ for $y_n:=x_n\circ\lambda_n$.
This also means $\|\Delta y_n-\Delta x\|\to 0$.

Fix a time $\tau \in [0,1]$ and assume  $\myAA_\tau^+(x)\neq \emptyset$.
Recall that $\myAA^+_\tau(x)$ is a finite set. 
Let $\myAA^+_\tau=\{s_1,\dots, s_N\}\neq \emptyset$  with $s_1= R_\tau(x)=\min\myAA^+_\tau(x)$.
Observe that $R_\tau(y_n)=\lambda_n^{-1}(R_\tau(x_n))$ and, thus,
since  $\|y_n-x\|\to 0$,
 $\CCCC(R_\tau(x_n))=\CCCC(R_\tau(y_n))$.
 For (ii) it thus suffices to show that $\CCCC(R_\tau(y_n))\subseteq \myAA^+_\tau(x)$.

To see this, note that there exist $\delta>0$ and  $n_0\in\NN$ such that for all $n\ge n_0$, $8(1+\#\myAA^+_\tau(x))\|y_n-x\|\le \delta$
(because  $\|y_n-x\|\to 0$ and  $\myAA^+_\tau(x)$ is finite),
and, also, 
\be\label{bound}
\SSS_{\Delta}\Big(x\!-\!\sum_{s\in\myAA^+_\tau(x)}\Delta x(s) {\bf 1}_{[s,1]}\Big)(\tau)< \Delta x(R_\tau(x))-\delta.
\ee
To explain \eqref{bound}: the quantity 
\ben
\Big(x\!-\!\sum_{s\in\myAA^+_\tau(x)}\Delta x(s)  {\bf 1}_{[s,1]}\Big)(\tau)
\een
is $x$ with all positive jumps equal in magnitude to the largest jump up till time $\tau$ subtracted.
Applying the operator $\SSS_{\Delta}$ to this produces the largest of the remaining jumps, hence the second largest jump, in magnitude, of $x$,
 up till time $\tau$.
This is strictly smaller than the magnitude of the largest jump, which is $ \Delta x(R_\tau(x))$. So indeed there is a $\delta>0$ such that \eqref{bound} holds.


Let $\alpha \in\CCCC(R_\tau(y_n))$ be the limit along a subsequence $(n')\subseteq(n)$.
Contrary to the hypothesis, suppose $\alpha \notin \myAA^+_\tau(x)$ and, thus, $\alpha_{n'}:=R_\tau(y_{n'})\notin \myAA^+_\tau(x)$ for all sufficiently large $n'$. For those $n'$, also being larger than $n_0$, observe that
\begin{eqnarray*}
\Delta x(R_\tau(x)) &=& \SSS_{\Delta}(x)(\tau)\cr
&\le&
 \SSS_{\Delta}(y_{n'})(\tau)+2\|y_{n'}-x\|\quad {\rm (by}\ \eqref{2in})  \\
&=&
\SSS_{\Delta}\Big(y_{n'}\!-\!\sum_{s\in\myAA^+_\tau(x)}\Delta y_{n'}(s) {\bf 1}_{[s,1]}\Big)(\tau)+2\|y_{n'}-x\|.
\eean
Here the second equality holds because $\alpha_{n'} \notin \myAA^+_\tau(x)$ 
implies ${\Delta}y_{n'}(s)\le  \SSS_{\Delta}(y_{n'})(\tau)$
for any $s\in  \myAA^+_\tau(x)$, for large $n'$, and thus subtracting any such jumps from $y_{n'}$ does not affect the value of $ \SSS_{\Delta}(y_{n'})$.
Using  \eqref{2in} again now gives 
\bean
\Delta x(R_\tau(x)) 
&\le &
\SSS_{\Delta}\big(x\!-\!\sum_{s\in\myAA^+_\tau(x)}\Delta x(s) {\bf 1}_{[s,1]}\big)(\tau)+4(1+\#\myAA^+_\tau (x))\|y_{n'}-x\|\cr
&\le&
\Delta x(R_\tau(x))-\delta + \frac 12\delta
= \Delta x(R_\tau(x))-\frac 12\delta,
\eean
where the last inequality holds  by \eqref{bound}.
This contradiction gives $\alpha\in\myAA^+_\tau(x)$, completing the proof that $\CCCC(R_\tau(y_n)) \subset \myAA^+_\tau(x)$.

(iii)\ For this,  suppose   again
$\|\lambda_n-I\|\vee \|y_n-x\|\to 0$ and, in addition,  $\#\myAA_\tau^+(x)=1$. We can take
$16\|y_n-x\|\le \delta$ and
\eqref{bound} now takes the form
\[
 \SSS_{\Delta}\big(x\!-\!\Delta x(R_\tau(x)) {\bf 1}_{[R_\tau(x),1]}\big)(\tau)<\Delta x(R_\tau(x)) - \delta,
\]
 for some  $\delta>0$ and all $n\ge n_0$.  
Suppose that $R_\tau(y_n)\neq R_\tau(x)$ for some $n$. Then we must have $n<n_0$, as otherwise
\bean
\Delta x(R_\tau(x)) 
&\le&
|\Delta y_n(R_\tau(x))|+2\|y_{n}-x\|
\le |\Delta y_n(R_\tau(y_n))|+2\|y_{n}-x\|\cr
&\le&
|\Delta x(R_\tau(y_n))|+4\|y_{n}-x\|
< \Delta x(R_\tau(x)) -\delta + \frac 13\delta\cr
&=& \Delta x(R_\tau(x)) -\frac 23\delta.
\eean
This  contradiction proves the result.     \halmos

\noindent {\bf Proof of Theorem \ref{signTies}:}\ 
Assume  $\sup_{\tau\in [0,1]}\#\BB_\tau(x) = 0$. We first show that $\TTTT_{rim}$ is $\Lambda$-compatible. 
Let $\lambda \in \Lambda$ and recall from Proposition \ref{operCon} that $\wt \SSS_\Delta$ is $\Lambda$-compatible.
Since
\[\wt L_\tau(x \circ \lambda) 
				   = \sup\{0\le s\le \tau: \, \Delta x(\lambda(s)) = \wt \SSS_\Delta(x) ( \lambda(\tau)) \} 
\] 
and 
\[ \wt L_{\lambda(\tau)}(x) =  \sup\{0\le s\le \lambda(\tau): \, \Delta x(s) = \wt \SSS_\Delta(x)  (\lambda(\tau)) \},
\] 
we have $\lambda^{-1} \circ \wt L_{\lambda(\tau)}(x) = \wt L_{\tau}(x\circ \lambda)$.
Thus
\begin{align*}
&
 \Delta(x\circ \lambda)(\wt L_{\tau}(x\circ \lambda))
 = \Delta (x\circ \lambda) (\lambda^{-1} \circ \wt L_{\lambda(\tau)} (x)) = \Delta x(\wt L_{\lambda(\tau)}(x)).
\end{align*}
Then $\TTTT_{rim}$ is $\Lambda$-compatible, because 
\bean
\TTTT_{rim}(x \circ \lambda) 
&=&
 x \circ \lambda -  \Delta(x\circ \lambda)(\wt L_{\tau}(x\circ \lambda))
= x \circ \lambda -  \Delta x(\wt L_\tau(x)) \circ \lambda \cr
&=&
 (x - \Delta x(\wt L_\tau(x)) \circ \lambda = \TTTT_{rim}(x) \circ \lambda.
\eean

It remains to show that  $\TTTT_{rim}$ is $\|\cdot\|$-continuous at $x$.
Suppose $\| x_n - x\| \to 0$. 
If $\wt \myAA_\tau(x) = \emptyset$, then $\TTTT_{rim}(x) = x$, hence $\TTTT_{rim}$ is trivially $\|\cdot\|$-continuous. Alternatively, suppose $\wt \myAA_\tau \neq \emptyset$.
Then 
 \[ \|\TTTT_{rim} (x_n) - \TTTT_{rim}(x)\|\le \| x_n-x\| + 
 \sup_{0\le \tau \le 1} \big| \Delta x_n(\wt L_\tau(x_n)) - \Delta x(\wt L_\tau(x))  \big|.
 \]
The first term on the RHS tends to $0$ and the second term on the RHS does not exceed
\begin{align}\label{ST1}
&\sup_{0\le \tau \le 1} \big| \Delta x_n(\wt L_\tau(x_n)) - \Delta x(\wt L_\tau(x_n))  \big|+ \sup_{0\le \tau \le 1} \big| \Delta x(\wt L_\tau(x_n)) - \Delta x(\wt L_\tau(x))  \big| \nonumber \\
 \le & 2 \|x_n -x\| + \sup_{0\le \tau \le 1} \big| \Delta x(\wt L_\tau(x_n)) - \Delta x(\wt L_\tau(x))  \big|. 
\end{align}
By Lemma \ref{lemaccpointsT}, $\CCCC_{\infty}(\wt L_\tau(x_n)) \subseteq \wt\myAA_\tau(x)$ for each $\tau \in [0,1]$. 
If $\wt L_{\tau}(x_n) \to \wt L_{\tau}(x)$, then the second term on the RHS of \eqref{ST1} tends to $0$. Suppose $\wt L_{\tau}(x_n) \to s_1 \neq \wt L_{\tau}(x) = s_2 $ where $s_1, s_2 \in \wt \myAA_\tau(x)$. Then $|\Delta x(s_1)| = |\Delta x(s_2)|$. 
But since $\# \BB_\tau(x) = 0$ for all $\tau \in [0,1]$, we also have $\Delta x(s_1) = \Delta x(s_2)$ for all $\tau \in [0,1]$. Thus $\sup_{0\le \tau \le 1} \big| \Delta x(\wt L_\tau(x_n)) - \Delta x(\wt L_\tau(x))  \big| \to 0$.  This completes the proof.
\halmos

\noindent {\bf Proof of Theorem \ref{theotrimcont}:}\ 
Again we only consider the case of $\myAA_1^+$.

\noindent 
{(i)}\
 Let $x_n\sko x$ in $\DD$ and  $\#\myAA_1^+(x)\le 1$.
 Then there exists a sequence $\lambda_n\in\Lambda$
such that $\|\lambda_n-I\|\vee \|x\circ \lambda_n-x\|\to 0$.
If $\myAA_1^+(x)=\emptyset$ then $\SSS_{\Delta} (x)\equiv 0$ and $ \RRR_{trim}(x)=x$, so
\begin{eqnarray*}
\| \RRR_{trim}(x_n)\circ\lambda_n- \RRR_{trim}(x)\|&=&\| \RRR_{trim}(x_n\circ\lambda_n)-x\|\\
&\le&\|x_n\circ\lambda_n-x\|+|\Delta (x_n\circ \lambda_n) (R_1 (x_n\circ \lambda_n))|\\
&= &\|x_n\circ\lambda_n-x\|+\|\SSS _{\Delta} (x_n\circ \lambda_n)-\SSS _{\Delta} (x)\|\\
&\le &3\|x_n\circ\lambda_n-x\|\to 0,\quad n\to\infty.
\end{eqnarray*}

Alternatively, if $\#\myAA_1^+ =1$, then by Lemma~\ref{lemaccpointsT} there is an $n_0\in\NN$ such
that $R_1(x_n\circ\lambda_n)=R_1(x)$ for all $n\ge n_0$ and, for those $n$, we also have 
\ben
 \RRR_{trim}(x_n)\circ \lambda_n= \RRR_{trim}(x_n\circ \lambda_n)=x_n\circ \lambda_n- \Delta x_n(R_1(x)){\bf 1}_{[R_1(x),1]}.
 \een
  As $n\to\infty$, the right hand-side converges uniformly (in the supremum norm) to $x-\Delta x(R_1(x)){\bf 1}_{[R_1(x),1]}= \RRR_{trim}(x)$, which 
shows that  $\RRR_{trim}$ is jointly $J_1$-continuous at $x$. 
For $r = 2$, recall the definition of $\RRR^{(2)}_{trim}$ in \eqref{defRt}. Since $\RRR_{trim}$ is $J_1$-continuous at $x$ and $\RRR_{trim}$ is assumed $J_1$-continuous at $\RRR_{trim}(x)$, then the composition 
$\RRR^{(2)}_{trim}(x) = \RRR_{trim}\big(\RRR_{trim}(x)\big)$ is $J_1$-continuous at $x$. An analogous argument holds for $r > 2$.

{(ii)} \
Contrary to the hypothesis, assume that $\{s_1,s_2\}\subseteq\myAA^+_1(x)$ for some $0\!<\!s_1\!<\!s_2\!\le\! 1$.
Noting that $s_1=R_1(x)$ and $ \RRR_{trim}(x)=x-\Delta x(s_1)\,{\bf 1}_{[s_1,1]}$ we introduce
\[x_n:=x+\frac{1}{n}\;{\bf 1}_{[s_2,1]},\qquad n\in\NN.\]
As $n\!\to\!\infty$, we have $\|x_n\!-\!x\|\!=\! 1/n\!\to\! 0$ and, in particular, $x_n\!\sko\! x$.
Observe that $R_1(x_n)=s_2$ and $ \RRR_{trim}(x_n)=x-\Delta x(s_2) {\bf 1}_{[s_2,1]}$, $n\in\NN$.
Hence, for all $n$,
\ben
\| \RRR_{trim}(x_n)- \RRR_{trim}(x)\|\ge | \RRR_{trim}(x_n)(s_1)- \RRR_{trim}(x)(s_1)|=|\Delta x(s_1)|>0.
\een

Finally, let $(\lambda_n)\subseteq\Lambda$ be such that $\lim_{n\to\infty}\|\lambda_n-I\|=0$. Then
$\lim_{n\to\infty}\lambda_n^{-1}(s_1)=s_1$. As $ \RRR_{trim}(x)$ is continuous at $s_1=R_1(x)$, 
 \[\delta_n\;:=\; \RRR_{trim}(x)(s_1)- \RRR_{trim}(x)(\lambda_n^{-1}(s_1))\to 0,\ {\rm as}\ n\to\infty,
 \]
and, thus,
\begin{eqnarray*}
\| \RRR_{trim}(x_n\circ \lambda_n)- \RRR_{trim}(x)\|&=&\| \RRR_{trim}(x_n)- \RRR_{trim}(x)\circ \lambda_n^{-1}\|\\
&\ge& | \RRR_{trim}(x_n)(s_1)- \RRR_{trim}(x) (\lambda_n^{-1}(s_1))|\\
&=& |\Delta x(s_1)+\delta_n|\to  |\Delta x(s_1)|>0,\ n\to\infty. \end{eqnarray*}
To summarise, we showed that $x_n\sko x$, but not $ \RRR_{trim}(x_n)\sko  \RRR_{trim}(x)$,  contradicting the $J_1$-continuity of $ \RRR_{trim}$ at~$x$.
\halmos

\section{Proof of Theorem \ref{maxtrim}}\label{s3}
Let $(\Delta X_s)_{0<s\le t}$ be the jumps of a L\'evy process $(X_s)_{0<s\le t}$ having L\'evy measure $\Pi_X$, with
ordered jumps $(\Delta X_t^{(i)})_{i\ge 1}$ and $(\Delta X_t^{(j),-})_{j\ge 1}$, as specified in Section \ref{s2}. 
 In what follows we will assume
 $\pibar_X^+(0+) =\pibar_X^-(0+)= \infty$ throughout, so there are always infinitely many positive and negative jumps of $X$, a.s., in any interval of time.

Let  $(\EEE_i)_{i\ge 1}$ be an i.i.d. sequence of exponentially distributed random variables with common parameter $E\EEE_i=1$ and
let $\Gamma_r:=\sum_{i=1}^r\EEE_i$ with $r\in\NN$. 
Write 
\ben  
\pibarpinv_X(x)=\inf\{y>0: \pibar^+_X (y) \le x\},\ x>0,
\een
for the right-continuous inverse of the right tail $\pibar^+_X$
(with similar notations for the left  tail $\pibar^-_X$ and the two sided 
 tail $\pibar_X$).
The following distributional equivalence can be deduced from Lemma 1.1 of BFM:
\be\label{ch}
\big(\Delta X_t^{(i)}\big)_{1\le i\le r}
\eqdr \big(\pibarpinv_X(\Gamma_i/t)\big)_{1\le i\le r},\ t>0,\ r\in\NN.
\ee
 We refer to BFM for  more background information on the properties of the extremal processes $(\Delta X_t^{(r)})_{t\ge 0}$ and the trimmed L\'evy processes.

 \noindent {\bf Proof of Theorem \ref{maxtrim}, Part (i)}\ 
 We give proofs just for 
$t\dto 0$; $t\to\infty$ is very similar.
Recall the definition of $\{\III_t(\tau)\}_{\tau \in [0,1]}$ in
\eqref{funcFC1} and assume the convergence of $\III_t$ to a stable process $\YYY$ as in \eqref{funcFC1}. The process $\YYY$ has L\'evy measure $\Pi_Y$ which is diffuse (continuous at each $x\in\R\setminus\{0\}$).


For each $\tau \in (0,1]$ the jump of $\III_t$ at $\tau$ is
\be\label{jumpi}
 \Delta \III_t(\tau) := \III_t (\tau) - \III_t(\tau-) = \frac{X_{t\tau} - \tau a_t}{b_t} -\frac{X_{t\tau-} - \tau a_t}{b_t}  = \frac{\Delta X_{t\tau}}{b_t}.
\ee
Hence, $\SSS_\Delta(\III_t)(\tau) = \SSS_\Delta (X_{t\tau}/b_t) $ for each $t > 0$, and we can write
\[ 
\Big(\frac{{}^{(r)}X^R_{\tau t} - \tau a_t}{b_t}\Big)_{\tau \in [0,1]}= \RRR^{(r)}_{trim}(\III_t).
\]
We want to apply the continuous mapping theorem and deduce the convergence in \eqref{ffunLook} from this.
By Theorem \ref{theotrimcont}, to apply the continuous mapping theorem it is enough to verify that there are no ties a.s. among the first $r$ largest positive jumps in the limit process $\YYY$. 
Let $\CCC :=\{ x \in \DD \,: \,  \#\myAA^+_1(\RRR_{trim}^{(j)}(x)) \le 1 \text{ for all } j = 0,\ldots, r-1 \}$. We wish to show that $\PP(\YYY \in \CCC) = 1$. 
Denote by $\Delta Y_1^{(j)}$ the $j$th largest jump of $\YYY$ up to time $1$.
Note that 
\[ \PP(\YYY \in \CCC) \ge 1 - \sum_{j=1}^{r}\PP(\Delta Y_1^{(j)} = \Delta Y_1^{(j+1)}).
\] 
Since  $\Pi_Y$ is diffuse, we have $\pibar^+_Y(\pibarpinv_Y(v)) = v = \pibar^+_Y(\pibarpinv_Y(v)-)$ for all $v >0$. Thus, by \eqref{ch} (with $X$ and $\Pi_X$ replaced by $Y$ and $\Pi_Y$), 
\begin{align*}
\PP(\Delta Y_1^{(j)} = \Delta Y_1^{(j+1)}) &= \PP\big\{ \pibarpinv_Y(\Gamma_j) = \pibarpinv_Y(\Gamma_j + \EEE_{j+1}) \big\} \nonumber \\
									&= \int_{0}^\infty \PP\big\{\pibarpinv_Y(v) = \pibarpinv_Y(v + \EEE_{j+1})\big\} \, e^{-v}\frac{v^{j-1}}{(j-1)!} \rmd v \nonumber \\
									&= \int_{0}^\infty \PP\big\{0 \le \EEE_{j+1} \le \pibar^+_Y(\pibarpinv_Y(v)-) -v\big\} \, e^{-v}\frac{v^{j-1}}{(j-1)!} \rmd v\cr
&= 0.
\end{align*}
So we can apply  Theorem \ref{theotrimcont} as forecast and complete the proof.
\halmos

 \noindent {\bf Proof of Theorem \ref{maxtrim}, Part (ii)}\  
We first prove \eqref{ffunLaw} and consider only the trimming operator $\TTT_{rim}^{(1,+)}$ ($\TTT_{rim}^{(r,s)}$ is treated analogously.)
By \eqref{jumpi} we can write, for each $\tau \in (0,1]$
and $t>0$,
\ben\label{n1}
\TTT_{rim}^{(1,+)}(\III_t)(\tau) = \frac{X_{\tau t} - \tau a_t}{b_t} - \sup_{0<s\le \tau} \Delta \III_t(s)= \frac{{}^{(1)}X_{\tau t} - \tau a_t}{b_t}. 
\een
 Since $\TTT_{rim}^{(1,+)}$ is globally $J_1$-continuous on $\DD$ by Proposition \ref{operCon}, we can apply the continuous mapping theorem to get  
that
\begin{equation*}\label{j1_4}
\left(\frac{{}^{(1)}X_{\tau t} - \tau a_t}{b_t}  \right)_{0\le \tau\le 1}=
 \TTT_{rim}^{(1,+)}(\III_t)  \todr  \TTT_{rim}^{(1,+)}(\YYY) =  
\big( {}^{(1)}Y_\tau\big)_{0\le \tau\le 1},
\end{equation*}
in $J_1$,  as $t \dto 0$ or $t\to\infty$.
This completes the proof of \eqref{ffunLaw}.

 \noindent {\bf Proof of Theorem \ref{maxtrim}, Part (iii)}\ 
Again by \eqref{jumpi}, we have  
\[   
\bigg(\frac{{}^{(r)}\wt X_{\tau t}-\tau a_t}{b_t}\bigg)_{ 0 \le \tau \le 1} =\TTTT^{(r)}_{rim}(\III_t).
\] 
Recall from Theorem \ref{signTies} that $\TTTT^{(r)}_{rim}$ is $J_1$-continuous on $\wt \CCC$, where
\[ 
\wt \CCC : = \{x \in  \DD \,:\, \# \BB_\tau(\TTTT^{(j)}_{rim}(x)) = 0 \,\, \text{for all } \tau \in [0,1], \, \, j = 0, \ldots, r-1 \}.
\]
Thus, in order to apply the continuous mapping theorem, we need to show that $\PP(\YYY \in \wt \CCC) = 1$.
Note that $\wt \CCC \supseteq \wt \VVV$, where 
\ben
\wt \VVV: = \{x \in \DD: \# \wt \myAA_\tau(\TTTT^{(j)}_{rim}(x)) \le 1 \ \text{for all } \tau \in [0,1],\,\, j = 0, \ldots, r-1 \}.
\een
Hence, it is enough to show that $\PP( \YYY \in \wt \VVV) = 1$, or,
equivalently,
\be\label{YW}
\PP\bigg(\bigcup_{1\le j \le r} \bigcup_{ 0< \tau < 1}\big\{|\wt{\Delta Y}_\tau^{(j+1)}|    = |\wt{\Delta Y}_\tau^{(j)}| \big\}\bigg) = 0,
\ee
where $\wt{\Delta Y}_\tau^{(j)}$ denotes the $j$th largest modulus jump of  $\YYY$ up till time $\tau$.

To simplify notation, during the remainder of this proof, write  $\Delta_t$ for the modulus jumps $|\Delta Y_t|$, and for their ordered values in the intervals $[0,t]$ or $[0,t)$, write $\Delta_t^{(j)}= |\wt{\Delta Y}_t^{(j)}|$ or  $\Delta_{t-}^{(j)}= |\wt{\Delta Y}_{t-}^{(j)}|$, $t>0$, $j=1,2,\ldots$. We aim to show
\be\label{Yv}
\PP\bigg(\bigcup_{ 0< \tau< 1}\{\Delta_\tau^{(j+1)}  = \Delta_\tau^{(j)} \}\bigg) = 0,\ j=1,2,\ldots,r,
\ee
from which \eqref{YW} will follow immediately. 

We consider first the case $j=1$.
Define a sequence of random times $(\tau_k)_{k\ge 0}$ by
\begin{equation}\label{taudef1}
\tau_0=1,\quad  {\rm and}\quad 
\tau_{k+1}:=\inf\{0<t<\tau_k: \Delta_t=\Delta_{\tau_k-}^{(1)}\},\ k=0,1,\ldots
\end{equation}
Since $\lim_{t\dto 0}\Delta_t^{(1)}=0$ a.s., we have
$0<\tau_{k+1}<\tau_k\le 1$ and $\lim_{k\to\infty}\tau_k=0$ a.s.
On  $\{\tau_{k+1}\le t<\tau_k\}$,  we have $\Delta_t^{(1)} = \Delta_{\tau_k-}^{(1)}$, hence on the event 
$\{ \Delta _t^{(2)}=\Delta_t^{(1)} \}$,
\[ 1=  \frac{\Delta_t^{(2)}}{ \Delta_t^{(1)}} \le \frac{\Delta_{\tau_k-}^{(2)}}{\Delta_{\tau_k-}^{(1)}} \le 1.
\]
This implies that 
\begin{align*}
& \bigcup_{ 0< \tau< 1} \{
\Delta _\tau^{(2)}=\Delta_\tau^{(1)}\}
= \bigcup_{k\ge 0}\bigcup_{ t \in [\tau_{k+1}, \tau_{k})}\{ \Delta _t^{(2)}=\Delta_t^{(1)} \} \nonumber \\
&=
 \bigcup_{k\ge 0} \{\Delta_{\tau_k-}^{(2)} 
 = \Delta_{\tau_k-}^{(1)} \} \nonumber \\
&=
\big\{\Delta_{t-}^{(2)} = \Delta_{t-}^{(1)}\ 
{\rm for\ some}\ t\le 1\ {\rm with}\ \Delta_t>  \Delta _{t-}^{(1)}\big\}=:E,\ {\rm say}.
\end{align*}
Define
$\mathbb{S} = \sum_{0 < t <1} \delta_{(t, \Delta_t)}$, where $\delta_{(t, \Delta_t)}$ denotes a point mass at $(t, \Delta_t)$.
$\mathbb{S}$ is a Poisson random measure on $(0,1)\times (0,\infty)$ with intensity $\rmd t \times \Pi_Y(\rmd x)$.
Let 
\ben\label{defN}
N=
\int_{(0,1)\times (0,\infty)}\textbf{1}\big\{\Delta ^{(2)}_{t-}=\Delta^{(1)}_{t-}< x\big\} \mathbb{S}(\rmd t\times \rmd x)
\een
be the number of points $(t,\Delta_t)$ which satisfy 
$\Delta^{(2)}_{t-}=\Delta^{(1)}_{t-}< 
\Delta_t$ with $ t< 1$. 
Then, recalling that  $\Delta_t^{(j)}= |\wt{\Delta Y}_t^{(j)}|$, event $E$ has probability 
\bea\label{A+B}
&&\PP(E)=
\PP(N\geq 1)
\leq \EE(N)\hspace{1cm}\mbox{(by Markov's inequality)}\cr
&=&
\int_0^1 \rmd t \int_{x>0}\EE \, \textbf{1}\big\{\Delta^{(2)}_{t-}= \Delta^{(1)}_{t-}< x\big\}\Pi_Y(\rmd x) \cr
&&\cr
&=&
\int_0^1\rmd t\int_{x>0}\int_{\pibarinv_Y(y/t)<x}
\PP\big(y+\EEE_{2}\le t\pibar_Y\big(\pibarinv_Y(y/t)-\big)\big)\PP\left(\EEE_1 \in \rmd y\right)\Pi_Y(\rmd x)
\cr
&&\cr
&=&
\int_0^1\rmd t\int_{x>0}\int_{\pibarinv_Y(y)<x}
\PP\left(\EEE_{2} \le t\left(\pibar_Y\left(\pibarinv_Y(y)-\right)-y\right)\right)e^{-ty}\rmd y\, \Pi_Y(\rmd x) \cr
&=&  0.
\eea
In the second equality we used the compensation formula,   and in the third we used a version of \eqref{ch} appropriate to the $ |\wt{\Delta Y}_t^{(j)}|$. The last expression in \eqref{A+B} is 0 because $\pibar_Y$ is diffuse, so  $\pibar_Y(\pibarinv_Y(v)) = v = \pibar_Y(\pibarinv_Y(v)-)$ for all $v >0$.
This means, with probability $1$, there are no tied values among the largest jumps in $(\Delta_\tau)_{0<\tau<t}$ for all $t \in (0,1)$.
(Note that this is ostensibly a much stronger statement than requiring there be no tied values among the largest jumps up until a fixed time $t$.)

Next we consider $j=2$. It is enough to show that $\PP\big(\bigcup_{ 0< t < 1} \{\Delta_t^{(3)}= \Delta _t^{(2)} \}\big) = 0$. We restrict ourselves to the event 
$\FFF : =\bigcap_{ 0< t < 1} \{\Delta _t^{(2)}\neq \Delta _t^{(1)} \}$, which we have proved has probability 1. On this event, there are no ties for the largest value among $(\Delta_t)_{0<t\le 1}$.
Recall the definition of the sequence $(\tau_k)$ in \eqref{taudef1}.  
The largest jump $\Delta_t^{(1)}$  remains constant on the interval $ \tau_{k+1}\le  t< \tau_k$. We aim to subdivide the interval $[\tau_{k+1}, \tau_k)$ so that the second largest jump up till time $t$, which is strictly less than $\Delta_{\tau_{k+1}}$, is constant within that subinterval. First we consider the case when $\Delta_t^{(2)} = \Delta_{\tau_{k+1}}$. Define for each $k \in \NN_0 := \NN \cup \{0\}$,  
\[ s_k := \sup\{ 0 < t < \tau_k :  \Delta_{\tau_{k+2}}= \Delta_t^{(2)}  \}.
\] Note that $s_k \ge \tau_{k+1}$ as $  \Delta_{\tau_{k+2} }= \Delta_{\tau_{k+1}}^{(2)}$. 
Next define a further sequence $(\sigma_{m}(k))_{m\ge 1}$ in $[s_k, \tau_{k})$ such that $\sigma_0(k) = \tau_{k}$ and for $m=1,2,\ldots,$
\[ 
\sigma_{m}(k) = \inf\{0< t < \sigma_{m-1}(k): \Delta_t = \Delta_{\sigma_{m-1}(k)-}^{(2)} \}\vee s_k.
\]
Then we can decompose
\begin{align}\label{iter0}
\bigcup_{ 0< t < 1} \{\Delta_t^{(3)}= \Delta _t^{(2)} \} 
= \bigcup_{k\ge 0}  \Big(\bigcup_{ \tau_{k+1}\le t < s_k }\,\, \medcup \,\, \bigcup_{m\ge 1} \bigcup_{ \sigma_{m}(k) \le t < \sigma_{m-1}(k) }  \Big)\{\Delta_t^{(3)} = \Delta_t^{(2)}\} .
\end{align}
When $\Delta_t^{(3)} = \Delta_t^{(2)}$ and $\{\tau_{k+1}\le t<s_k\}$,  we have $\Delta_t^{(2)} = \Delta_{s_k-}^{(2)} = \Delta_{\tau_{k+2}}$, hence 
\[ 1=  \frac{\Delta_t^{(3)}}{ \Delta_t^{(2)}} \le \frac{\Delta_{s_k-}^{(3)}}{\Delta_{s_k-}^{(2)}} \le 1.
\]
When $\Delta_t^{(3)} = \Delta_t^{(2)}$ and $\{\sigma_{m+1}(k)\le t<\sigma_{m}(k)\}$, $m\ge 1, k\in \NN_0$, we have $\Delta_t^{(2)} = \Delta_{\sigma_m(k)-}^{(2)} $, hence 
\[ 1=  \frac{\Delta_t^{(3)}}{ \Delta_t^{(2)}} \le \frac{\Delta_{\sigma_{m}(k)-}^{(3)}}{\Delta_{\sigma_{m}(k)-}^{(2)}} \le 1.
\]

So the events on the RHS of \eqref{iter0} are subsets of 
\begin{align}\label{iter1}
 & \bigcup_{k\ge 0}  \big\{ \Delta_{t-}^{(3)} = \Delta _{t-}^{(2)} \,\,  {\rm  for\ some}\ t \in [\tau_{k+1}, \tau_k)  \ {\rm with}\ \Delta_t> \Delta _{t-}^{(2)}\big\} \nonumber \\
= & \bigcup_{k\ge 0}\bigcup_{t \in [\tau_{k+1}, \tau_k)} \big\{ \Delta_{t-}^{(3)}= \Delta_{t-}^{(2)}, \, \Delta_t > \Delta_{t-}^{(2)} \big\} \nonumber \\
=& \big \{\Delta_{t-}^{(3)}= \Delta_{t-}^{(2)} \text{ for some } t < 1 \, \text{with } \Delta_t > \Delta_{t-}^{(2)} \big\}.
\end{align}
The probability of the event on the RHS of \eqref{iter1} can be computed in a similar way as in \eqref{A+B} to be 0.
Hence, reverting to the original notation,  we have
\[
 \PP\Big( \bigcup_{ 0< t < 1} \{|\wt{\Delta Y}_t^{(3)}|  = |\wt{\Delta Y}_t^{(2)}| \}\Big) 
\le  \, \PP\Big( \bigcup_{ 0< t < 1} \{|\wt{\Delta Y}_t^{(3)}|  = |\wt{\Delta Y}_t^{(2)}| \}, \, \FFF \Big)  +  \PP(\FFF^c) =0. 
\]
For $j \ge 3$, we can proceed iteratively with similar arguments to arrive at \eqref{Yv} hence \eqref{YW}.
This completes the proof of \eqref{ffunLaw2}.
\halmos


\section{Applications to Reinsurance}\label{S6}

Many examples can be generated from the convergences in
\eqref{ffunLook}, \eqref{ffunLaw} and \eqref{ffunLaw2} using  the continuous mapping theorem.
Here we mention one which is of particular interest in reinsurance.
The idea is that the largest claim up to a specified time incurred by an insurance company (the ``cedant") is referred to a higher level insurer (the ``reinsurer").  See Fan et al. \cite{FGMSW2016} for details and references to the applications literature. This is known as the  largest claim reinsurance (LCR) treaty: having set a fixed follow-up time $t$, we delete from the process the largest claim occurring up to and including that time. We refer to Ladoucette and Teugels \cite{teugels2006} and Teugels \cite{teugels2003} for more detailed expositions.

The LCR procedure can be made prospective by implementing it as a forward looking dynamic procedure in real time, from the cedant's point of view. Designate as time zero the time at which the reinsurance is taken out. The first claim arriving after time $0$ is referred to the reinsurer and not debited to the cedant. Subsequent claims smaller than the initial claim are paid by the cedant until a claim larger than the first (the previous largest) arrives. The difference between these two claims is referred to the reinsurer and not debited to the cedant.  The process continues in this way so that at time $t$, the accumulated amount referred to the reinsurer equals the largest claim up till that time. This procedure has the same effect as applying the ``trim as you go'' operator to the risk 
process. (It is also possible to apply ``lookback" trimming to a reinsurance model in a natural way.)

A primary quantity of interest is the ruin time, at which the process $(X_t)_{t\ge 0}$ describing the claims incoming to the company reaches a high level, $u>0$. After reinsurance of the $r$ highest claims, the process is reduced to $({}^{(r)}X_t)_{t\ge 0}$, with ruin time
$T^{(r)}(u)= \inf\{t>0: {}^{(r)}X_t>u\}$.
Supposing $X$ is L\'evy with heavy tailed canonical measure, $\Pi_X$, not uncommon assumptions in the modern insurance literature, we assume 
\eqref{dom} with no centering necessary, and from the continuous mapping
 theorem immediately deduce an asymptotic distribution for 
 $\sup_{0<\tau\le 1}{}^{(r)}X_{t\tau}/b_t$ as $t\to\infty$,   and hence for $T^{(r)}(\cdot)$, for high levels.  Specifically,
if  \eqref{ffunLaw} holds with $t\to\infty$ and $s=0$, then
\ben
\lim_{t\to\infty}\PP\left(T^{(r)}(ub_t)>t\right)
= \PP\left( T_Y^{(r)}>u\right),\ u>0,
\een
where $T_Y^{(r)}= \inf\{t>0: {}^{(r)}Y_t>1\}$.

\medskip\noindent {\bf Acknowledgement}\
We are grateful to Prof. Jean Jacod whose suggestion to YFI
 at the Building Bridges conference in Braunschweig 2013 started us on this pathwise analysis.
Buchmann and Maller's  research was partially funded by ARC Grants DP1092502 and  DP160104737.  Yuguang Ipsen was formerly at the  Australian Research Council Centre of Excellence for Mathematical and Statistical Frontiers, School of Mathematics and Statistics, University of Melbourne.
She acknowledges support from the Australian Research Council.
We are grateful also to two referees whose careful readings resulted in significant improvements to the paper.

{\footnotesize \bibliography{Library_Levy_Mar2017}}
\bibliographystyle{apt}

\end{document}